\documentstyle[12pt,amscd]{amsart}
\newtheorem{lm}{Lemma}

\newcommand{\proj}{\bold P}

\newcommand{\barr}{\overline}

\newcommand{\rarr}{\rightarrow}
\newcommand{\oh}{{\cal{O}}}
\newcommand{\M}{\barr{M}_C(d)}
\newcommand{\Me}{\barr{M}_1(d)}
\newcommand{\eqq}{\stackrel{\sim}{=}}
\newcommand{\deli}{\bigtriangleup}
\newcommand{\com}{\Bbb{C}}
\begin{document}
\title{A Note On Elliptic Plane Curves With Fixed $j$-Invariant}
\author{Rahul Pandharipande}
\date{10 May 1995}
\maketitle
\pagestyle{plain}
\section{Summary}
Let $N_d$ be the number of irreducible, reduced, nodal, degree
$d$ {\em rational plane curves} passing through $3d-1$ general points in the
complex projective plane $\proj^2$. The numbers $N_d$ satisfy a
beautiful recursion relation ([K-M], [R-T]):
$$N_1=1$$
$$\forall d>1, \ \ \ N_d= \sum_{i+j=d, \ i,j>0}
N_iN_j \bigg( i^2j^2 {3d-4 \choose 3i-2} - i^3j {3d-4 \choose 3i-1} \bigg).$$
Let $E_{d,j}$ be the number of irreducible, reduced, nodal,
degree $d$, {\em elliptic plane curves with fixed j-invariant j}
passing through $3d-1$ general points $\proj^2$.
$E_{d,j}$ is defined for $d\geq 3$ and $\infty \neq j \in \barr{M}_{1,1}$.
In this note, the following relations are established:
\begin{eqnarray*}
\forall j\neq 0,1728,\infty,  &  E_{d,j}= {d-1\choose 2}N_d, \\
j=0, & E_{d,0}= {1\over 3} {d-1 \choose 2} N_d, \\
j=1728, & \ E_{d,1728}= {1\over 2} {d-1\choose 2} N_d.
\end{eqnarray*}
If $d\equiv 0 \ mod \ 3$, then $ 3 \not\mid {d-1 \choose 2}$.
Since $E_{3\hat{d},0}$ is an integer,
$N_{3\hat{d}} \equiv 0 \ mod \ 3$ for $\hat{d}\geq 1$.
In fact, a check of values in [D-I] shows
$N_d \equiv 0 \ mod \ 3$ if and only if $d \equiv 0 \ mod \ 3$
for $3\leq d \leq 12$. P. Aluffi has calculated $E_{3,j}$ for
$j<\infty$ in [A]. Aluffi's results agree with the above formulas.

Thanks are due to Y. Ruan for discussions on
Gromov-Witten invariants and quantum cohomology. The question
of determining the numbers $E_{d,j}$ was first considered
by the author in a conversation with Y. Ruan.

\section{Kontsevich's Space of Stable Maps}
\subsection{ The Quasi-Projective Subvarieties $U_C(\Gamma,c,\barr{w})$
, $U_{j=\infty}(\Lambda,\barr{w})$}
\label{qpd} Fix $d\geq 3$ for the entire paper.
Let $C$ be a nonsingular elliptic curve or an irreducible, nodal
rational curve of arithmetic genus 1.
Consider the coarse moduli space of $3d-1$-pointed stable maps
from $C$ to $\proj^2$ of degree $d\geq 3$,
$\barr{M}_{C, 3d-1}(\proj^2, d)$. For convenience, the
notation $\M = \barr{M}_{C,3d-1}(\proj^2,d)$ will be used.
Let $S_d= \{1,2,\ldots,3d-1\}$ be the marking set. Constructions of
$\M$ can be found in [Al], [K], [P].
Let $\Gamma$ be a tree consisting of
a distinguished vertex $c$, $k\geq 0$ {\em other} vertices
$v_1, \ldots, v_k$, and $3d-1$ marked legs.
Let $0\leq e \leq d$. Weight the vertex $c$ by $e$.
Let $w_1, \ldots, w_k$ be non-negative integral weights
of the vertices $v_1, \ldots, v_k$ satisfying
$$e+ w_1 + \cdots + w_k =d.$$
Denote the weighting by $\barr{w}=(e,w_1,\ldots, w_k)$.
The marked, weighted tree with distinguished vertex
$(\Gamma, c, \barr{w})$ is {stable} if the following
implication holds for all $1\leq i\leq k$:
$$w_i=0 \ \ \Rightarrow \ \ valence(v_i)\geq 3.$$
Two marked, weighted trees
with distinguished vertex $(\Gamma, c,\barr{w})$ and $(\Gamma', c',\barr{w}')$
are isomorphic if there is an isomorphism of marked trees $\Gamma\rarr \Gamma'$
sending $c$ to $c'$ and respecting the weights.

A quasi-projective subvariety $U_C(\Gamma,c, \barr{w})$ of $\M$
is associated to each isomorphism class
of stable, marked, weighted graph with distinguished vertex $(\Gamma,
c,\barr{w})$.
$U_C(\Gamma,c,\barr{w})$ consists of stable maps
$\mu: (D, p_1,\ldots, p_{3d-1})
\rarr \proj^2$ satisfying the following conditions. The domain $D$ is
equal to a union:
$$D=C \cup \proj^1_1 \cup \cdots \cup \proj^1_k.$$ The
marked, weighted dual graph with distinguished vertex of the map $\mu$ is
isomorphic to $(\Gamma, c, \barr{w})$. The distinguished vertex of the
dual graph of $\mu$ corresponds to the (unique) component of $D$ isomorphic to
$C$.
Weights of the
dual graph of $\mu$ are obtained by the degree of $\mu$ on the components.
Note $U_C(\Gamma,c,\barr{w})=\emptyset$ if and only if $e=1$.

Let $(\Gamma,c,\barr{w})$ be a stable, marked, weighted tree
with distinguished vertex. Assume $e\neq 1$. The dimension of
$U_{C}(\Gamma,c,\barr{w})$ is determined as follows.
If $e\geq 2$, then
$$dim \ U_{C}(\Gamma,c,\barr{w})= 6d-2-k.$$
If $e=0$, then
$$dim \ U_{C}(\Gamma,c,\barr{w})= 6d-k$$
(where $k$ is the number of non-distinguished vertices of $\Gamma$).
These calculations are straightforward.

Let $C$ be a nonsingular elliptic curve. Every stable
map in $\M$
has domain obtained by attaching a finite number of marked trees to $C$.
By the definition of tree and map stability:
$$\bigcup_{(\Gamma,c, \barr{w})} U_C(\Gamma,c, \barr{w}) \ = \ \M.$$

Let $C$ be an irreducible, 1-nodal rational curve.
The quasi-projective varieties $U_C(\Gamma,c,\barr{w})$
do not cover $\M$. The curve $C$ can degenerate
into a simple circuit of $\proj^1$'s.
Let $\Lambda$ be a graph with 1 circuit ($1^{st}$ Betti number equal to $1$,
no self edges),
$k\geq 1$ vertices $v_0,\ldots, v_k$, and $3d-1$ marked legs.
Note the different vertex numbering convention. At least $2$ vertices
are required to make a circuit, so $k\geq 1$.
Let $w_0, w_1,\ldots,w_k$ be non-negative, integral weights summing to $d$.
The marked, weighted graph with $1$ circuit $(\Lambda, \barr{w})$ is
stable if each zero weighted vertex has valence at least 3.
A quasi-projective subvariety $U_C(\Lambda, \barr{w})$
of $\M$ is associated to each
isomorphism class of stable, marked, weighted graph with $1$ circuit
$(\Lambda, \barr{w})$. $U_C(\Lambda, \barr{w})$ consists of
stable maps with marked, weighted dual graphs isomorphic to
$(\Lambda, \barr{w})$. The union
$$\bigcup_{(\Gamma,c, \barr{w})} U_C(\Gamma,c, \barr{w}) \
\ \cup \ \ \bigcup_{(\Lambda, \barr{w})} U_C(\Lambda, \barr{w})\  \ = \ \M$$
holds by the definition of stability.

Finally, the dimensions of the loci $U_C(\Lambda,\barr{w})$
will be required. Let $(\Lambda, \barr{w})$ be a stable, marked,
weighted graph with $1$ circuit. Let $c_1,\ldots, c_l$ be the
unique circuit of vertices of $\Lambda$.
Let $e$ be the sum of the weights of the circuit vertices.
$U_C(\Lambda, \barr{w})=\emptyset$ if and only if
$e=1$.
If $e\geq 2$, then
$$dim\  U_C(\Lambda,\barr{w})= 6d-2-k.$$
If $e=0$, then
$$dim\ U_C(\Lambda, \barr{w})=  6d-k$$
(where $k+1$ is the total number of vertices of $\Lambda$).
Again, these results are straightforward.

\subsection{The Component $\barr{W}_{1}(d)$}
Let $\Me = \barr{M}_{1,3d-1}(\proj^2,d)$ be
Kontsevich's space of $3d-1$-pointed stable maps
from genus $1$ curves to $\proj^2$. There is canonical morphism
$$\pi: \Me \rarr \barr{M}_{1,1}$$
obtained by forgetting the map and all the markings except $1\in S_d$ (the
$3d-1$ possible choices of marking in $S_d$ all yield the
same morphism $\pi$).
Let $j\in \barr{M}_{1,1}$. By the
universal properties of the moduli spaces, there is a canonical bijection
$$\barr{M}_{C_j}(d) \rarr \pi^{-1}(j)$$
where $C_j$ is the elliptic curve (possibly nodal
rational) with $j$-invariant $j$.
When $j\in \barr{M}_{1,1}$ is automorphism-free, this bijection
is an isomorphism. For $j=0,1728$, the scheme theoretic
fiber $\pi^{-1}(j)$ is nonreduced.
Define an open locus $W_{1}(d) \subset \Me$ by
$[\mu: (D, p_1, \ldots, p_{3d-1}) \rarr \proj^2]\in W_{1}(d)$ if and only if
$D$ is irreducible.
By considering the natural tautological spaces over the
universal Picard variety of degree $d$ {\em line bundles} over
$M_{1,3d-1}$, it is easily seen that $W_{1}(d)$ is a
reduced, irreducible open set of dimension
$6d-1$. Let $\barr{W}_1(d)$ be the
closure of $W_{1}(d)$ in $\Me$.

\section{A Deformation Result}
Let $\Phi$ be a stable, marked, weighted tree with
distinguished vertex $c$ determined by the data: $k=1$, $(e,w_1)=(0,d)$.
There are $2^{3d-1}$ isomorphism classes of such $\Phi$
determined by the marking distribution.
Let $j\in \barr{M}_{1,1}$. The dimension of
$U_j(\Phi)$ is $6d-1$.
A point $[\mu]\in U_j(\Phi)$ has domain $C_j \cup \proj^1$.
There are $3d-1$ dimensions of the map $\mu|_{\proj^1}:\proj^1\rarr \proj^2$.
The incidence point $p= C_j\cap \proj^1$ moves in a $1$-dimensional family on
$\proj^1$. The remaining $3d-1$ markings move in $3d-1$ dimensions on $C_j$ and
$\proj^1$ (specified by the marking distribution).
$6d-1=3d-1+1+3d-1$.
A technical result is
needed in the computation of the numbers $E_{d,j}$.
\begin{lm}
\label{aa}
Let $I(\Phi,j)= \barr{W}_1(d) \cap U_j(\Phi) \subset \Me$.
The dimension of $I(\Phi,j)$ is bounded by $dim\ I(\Phi,j) \leq 6d-3$.
\end{lm}
\begin{pf}
Let $[\mu]\in I(\Phi,j)$ be a point. Let $D=C_j \cup \proj^1$ be
the domain of $\mu$ as above. The following condition will be shown to hold:
the linear series on $\proj^1$ determined
by $\mu|_{\proj^1}$ has vanishing sequence $\{0, \geq 2,*\}$ at the
incidence point $p=C_j \cap \proj^1$. The existence of
a point with vanishing sequence $\{0,\geq 2,*\}$ is a $1$-dimensional
condition on the linear series. The condition that the incidence
point $p$ has this vanishing sequence is an additional
$1$-dimensional constraint on $p$. Therefore, the dimension of $I(\Phi,j)$
is at most $6d-1-1-1=6d-3$. The vanishing sequence $\{0, \geq 2, *\}$
is equivalent to $d(\mu|_{\proj^1})=0$ at $p$.

It remains to establish the vanishing sequence $\{0,\geq 2,*\}$ at
$p$. This result is easily seen in explicit holomorphic coordinates.
Let $\deli_t$ be a disk at the origin in $\com$ with coordinate $t$.
Let $\eta: \cal{E} \rarr \deli_t$ be a flat family of
curves of arithmetic genus $1$ satisfying:
\begin{enumerate}
\item[(i.)] $\eta^{-1}(0)\eqq C_j$.
\item[(ii.)] $\eta^{-1}(t\neq 0)$ is irreducible, reduced, and (at worst)
nodal.
\end{enumerate}
For each $1\leq i \leq d$, let $\cal{G}_i=\cal{H}_i\subset \cal{E}$ be the
open subset of $\cal{E}$ on which the morphism $\eta$ is {\em smooth}.
Consider the fiber product:
$$X= \cal{G}_1 \times_{\deli_t} \cdots \times_{\deli_t} \cal{G}_d
\times_{\deli_t} \cal{H}_1 \times_{\deli_t} \cdots \times_{\deli_t}
\cal{H}_d.$$
$X$ is a nonsingular open set of the $2d$-fold fiber product
of $\cal{E}$ over $\deli_t$. Let $Y\subset X$ be the
subset of points $y=(g_1,\ldots, g_d,h_1, \ldots, h_d)$ where
the two divisors $\sum g_i$ and $\sum h_i$ are linearly
equivalent on the curve $\eta^{-1}(\eta(y))$.
$Y$ is a nonsingular divisor in $X$.

Let $p\in C_j=\eta^{-1}(0)$ be a nonsingular point of $C_j$.
Certainly $p\in \cal{G}_i, \cal{H}_i$ for all $i$.
Let $\gamma: \deli_t \rarr \cal{E}$ be any local
holomorphic section of $\eta$ such that $\gamma(0)=p$.
Let $V$ be a local holomorphic field of vertical tangent vectors
to $\cal{E}$ on an open set containing $p$. The section $\gamma$ and the
vertical vector field $V$ together determine
local holomorphic coordinates $(t,v)$ on $\cal{E}$ at $p$.
Let $\phi_{V}:\cal{E}\times \com \rarr \cal{E}$ be the
holomorphic flow of $V$ defined locally near $(p,0)\in \cal{E}\times \com$.
The coordinate map $$\psi:(t,v) \rarr \cal{E}$$
is determined by $\psi(t,v)= \phi_V(\gamma(t),v)$.

Local coordinates on $X$ near the point $x_p=(p,\ldots, p,p,\ldots, p)\in X$
are given by $$(t,v_1,\ldots, v_d, w_1, \ldots w_d).$$
The coordinate map is determined by:
$$\psi_X(t,v_1,\ldots, v_d, w_1, \ldots w_d)=
(\psi(t,v_1),\ldots, \psi(t,v_d), \psi(t,w_1), \ldots, \psi(t,w_d))\in X.$$
Note $x_p\in Y$. Let $f(t,v_1,\ldots, v_d, w_1, \ldots w_d)$ be a
local equation of $Y$ at $x_p$.
Since $f$ is identically $0$ on the line $(t,0,\ldots,0,0,\ldots,0)$,
\begin{equation}
\label{tvan}
\forall k\geq 0, \ \ {\partial ^k f\over \partial t^k} |_{x_p} =0.
\end{equation}
The tangent directions in the plane $t=0$ correspond to
divisors on the fixed curve $C_j$. Here, it is well know (up to $\com^*$-
factor)
\begin{equation}
\label{lion}
{\partial f\over \partial v_i}|_{x_p}=+1, \ \
   {\partial f\over \partial w_i}|_{x_p}=-1.
\end{equation}
Equations (\ref{tvan}) and (\ref{lion})
are the only properties of $f$ that will be used.

Let $\hat{\eta}:\cal{\hat{E}} \rarr \deli_t$ be the family obtained
by blowing-up $\cal{E}$ at $p$ and adding $3d-1$-marking.
Let $\mu:\cal{\hat{E}}\rarr \proj^2$ be a morphism.
Let $\hat{\eta}^{-1}(0)=D= C_j \cup \proj^1$.
Assume the following conditions are satisfied:
\begin{enumerate}
\item[(i.)] $\mu$, $\hat{\eta}$, and the $3d-1$ markings determine
a family of Kontsevich stable pointed maps to $\proj^2$.
\item[(ii.)] The markings of $D$ are distributed according to $\Phi$.
\item[(iii.)] $deg(\mu|_{C_j})=0$, $deg(\mu|_{\proj^1})=d$.
\end{enumerate}
Let $L_1$, $L_2$ be general divisors  of $\mu^*(\oh_{\proj^2}(1))$
that each intersect $\proj^1$ transversely at $d$ distinct points.
For $1\leq \alpha \leq 2$,
$L_{\alpha}$ breaks into holomorphic sections
$s_{\alpha,1}+\ldots+s_{\alpha,d}$ of $\hat{\eta}$ over a holomorphic disk
at $0\in \deli_t$. These sections $s_{\alpha,i}$ ($1\leq \alpha \leq 2$,
$1\leq i \leq d$) determine a map $\lambda: \deli_t \rarr Y$ locally at $0\in
\deli_t$.
Let an affine coordinate on $\proj^1$ be given by $\xi$ corresponding
to the normal direction
\begin{equation}
\label{cord}
{d\gamma\over dt}|_{t=0} + \xi \cdot V(p).
\end{equation}
Let $s_{1,i}(0)=\nu_i \in \com \subset \proj^1$,
$s_{2,i}(0)=\omega_i \in \com \subset \proj^1$ be given
in terms of the affine coordinate $\xi$.
The map $\lambda$ has the form
$$\lambda(t)= (t, \nu_1 t,\ldots, \nu_d t, \omega_1 t, \ldots, \omega_d t)$$
to first order in $t$ (written in the coordinates determined by $\psi_X$).
Equations (\ref{tvan}), (\ref{lion}), and the condition $f(\lambda(t))=0$
implies
\begin{equation}
\label{cony}
\sum_{i=1}^{d} \nu_i = \sum_{1}^{d} \omega_i.
\end{equation}
$L_1 \cap \proj^1$ is a degree $d$ polynomial with roots at
$\nu_i$. Condition (\ref{cony}) implies that the
sums of the roots (in the coordinates (\ref{cord}))
of general elements of the linear
series $\mu|_{\proj^1}$ are the same. Therefore, a constant
$K$ exists with the following property. If
$$\beta_0 + \beta_1 \xi+ \ldots + \beta_{d-1} \xi^{d-1} +\beta_{d} \xi^{d}$$
is an element of the linear series $\mu|_{\proj^1}$, then
$\beta_{d-1}+ K\cdot \beta_d=0$. The vanishing sequence at $\xi=\infty$
is therefore $\{0, \geq 2,*\}$. The point $\xi=\infty$ is the intersection
$C_j \cap \proj^1$.

Suppose $\tilde{\eta}:\cal{\tilde{E}} \rarr \deli_t$ is obtained
from $\cal{E}$ by a sequence of $n$ blow-ups over $p$. The fiber
$\tilde{\eta}^{-1}(0)$ is assumed to be $C_j$ union
a chain of $\proj^1$'s of length $n$. Each blow-up occurs in the
exceptional divisor of the previous blow-up. Let $\proj$ denote the
extreme exceptional divisor. Let $\mu: \cal{\tilde{E}} \rarr \proj^2$
be of degree $d$ on $\proj$ and degree $0$ on the other components
of the special fiber $\tilde{\eta}^{-1}(0)$. Let there be $3d-1$ markings
as before. It must be again concluded that the linear series on
$\proj$ has vanishing sequence $\{0, \geq 2,*\}$ at the node.

Let $\gamma$ be section of $\eta$ such that the lift of $\gamma$ to
$\tilde{\eta}$ meets $P$. Let the coordinates $(t,v)$ on $\cal{E}$
be determined by this $\gamma$ (and any $V$). An affine coordinate
$\xi$ is obtained on $\proj$ in the follow manner.
Let $\gamma_{\xi}$ be
the section of $\eta$ determined in $(t,v)$ coordinates by
$$\gamma_{\xi}(t)=(t,\xi t^n).$$
Let $\tilde{\gamma}_{\xi}$ be the lift of $\gamma_{\xi}$ to
a section of $\tilde{\eta}$.
The association $$\com \ni \xi \mapsto \tilde{\eta}(0)\in \proj$$
is an affine coordinate on $\proj$.
Let $L_1, L_2$ be divisors in the linear series
$\mu$ intersecting $\proj$ transversely. As before,
$L_{\alpha}$ breaks into
holomorphic sections $s_{\alpha,1}$. Let
$s_{1,i}=\nu_i \in \com\subset \proj$, $s_{2,i}=\omega_i \in \com \subset
\proj$.
As before, a map $\lambda:\deli_t \rarr Y$ is obtained from the
sections $s_{\alpha,i}$. In the coordinates determined by $\psi_X$,
$$\lambda(t)=(t, \nu_1 t^n+O(t^{n+1}), \ldots, \nu_d t^n+ O(t^{n+1}),
 \omega_1 t^n+O(t^{n+1}), \ldots, \omega_d t^n+O(t^{n+1})).$$
As before $f(\lambda(t))=0$. The term of leading order in $t$ of
$f(\lambda(t))$ is
$$ (\sum_{i=1}^{d} \nu_i- \sum_{i=1}^{d}\omega_i) \cdot t^n.$$
This follows from equations (\ref{tvan}) and (\ref{lion}).
The vanishing sequence $\{0,\geq 2,*\}$ is obtained as before.

By definition,
an element $[\mu]\in I(\Phi,j)$ can be obtained as the special fiber
of family of Kontsevich stable maps where
the domain is a smoothing of the node $p$.
After resolving the singularity in the total space at the node
$p$ by blowing-up, a family $\cal{\tilde{E}}$ is obtained.
The above results show the linear series on $\proj^1$ has
vanishing sequence $\{0,\geq 2,*\}$ at $p$.
\end{pf}

The markings play no role in the preceding proof. An identical
argument establishes the following:
\begin{lm}
\label{bb}
Let $\Phi$ be a stable, marked, weighted
tree with distinguished vertex satisfying $e=0$ and $w_i=d$ for some $i$.
Let $k$ be the number of non-distinguished vertices of $\Phi$.
Let $j\in \barr{M}_{1,1}$. Let $I(\Phi,j)= \barr{W}_1(d) \cap U_j(\Phi)$.
The dimension of $I(\Phi,j)$ is bounded by
$dim\ I(\Phi,j)\leq 6d-k-2$
\end{lm}
\begin{lm}
\label{cc}
Let $\Omega$ be a stable, marked, weighted graph with $1$ circuit.
Let $v_i$ be a non-circuit vertex with
weight $w_i=d$ (this implies $e=0$).
Let $k+1$ be the total number of vertices of $\Omega$.
Let $I(\Omega,\infty)=\barr{W}_1(d) \cap U_{\infty}(\Omega)$.
The dimension of $I(\Omega,\infty)$ is bounded by
$dim \ I(\Omega,\infty)\leq 6d-k-2$.
\end{lm}
\noindent
The vanishing sequence $\{0,\geq 2, *\}$ condition reduces
the dimensions of $U_C(\Phi)$, $U_{\infty}(\Omega)$ by $2$.

\section{The Numbers $E_{d,j}$}
The space of maps $\Me$ is equipped with $3d-1$ evaluation
maps corresponding to the marked points.
For $i\in S_d$, let $e_i:\barr{W}_1(d)\rarr \proj^2$ be the
restriction of the $i^{th}$ evaluation map to $\barr{W}_1(d)$.
let $\cal{L}_i= e_i^*(\oh_{\proj^2})$
Let
$$Z=c_1(\cal{L}_1)^2 \cap \ldots \cap c_1(\cal{L}_{3d-1})^2$$
Let $\pi_{\barr{W}}:
\barr{W}_{1}(d) \rarr \barr{M}_{1,1}\cong \proj^1$ be the
restriction of $\pi$ to $\barr{W}_1(d)$.
Let $$T=c_1(\pi_{\barr{W}}^*(\oh_{\proj^1}(1))).$$
Note $\barr{W}_{1}(d)$ is an irreducible, projective scheme of
dimension $6d-1$. The top intersection of line bundles on $\barr{W}_1(d)$,
$Z\cap T$,
is an integer.
\begin{lm}
\label{pal}
\begin{eqnarray*}
\forall j\neq 0,1728,\infty, & Z\cap T=E_{d,j}\ , \\
j=0, & \ Z\cap T = 3\cdot E_{d,0}\ , \\
j=1728, & \ \ Z\cap T= 2\cdot E_{d,1728}\ .
\end{eqnarray*}
\end{lm}
\begin{pf}
Via pull-back, lines in $\proj^2$ yield representative classes of
$c_1(\cal{L}_i)$. Therefore $3d-1$ general points in $\proj^2$,
$\barr{x}=(x_1, \ldots, x_{3d-1})$,
determine a representative cycle $Z_{\barr{x}}$ of the the class $Z$.
Let $\infty > j \in \barr{M}_{1,1}$.
Let $\pi_W$ be the restriction of $\pi$ to
$W_1(d)$.
It is first established for a general representative $Z_{\barr{x}}$,
\begin{equation}
\label{erst}
Z_{\barr{x}} \cap \pi_{\barr{W}}^{-1}(j) \subset \pi_W^{-1}(j).
\end{equation}
The statement (\ref{erst}) is proven by considering
the quasi-projective strata of $\barr{M}_{C_j}(d)$.

Note $\pi_W^{-1}(j)$ is the strata $U_{C_j}(\Gamma,c,\barr{w})$ where
$(\Gamma,c,\barr{w})$ is the trivial, stable, marked, weighted tree
with distinguished vertex. Assume now $(\Gamma,c,\barr{w})$ is
not the trivial tree.
By the equations for the dimension of $(\Gamma,c,\barr{w})$
$$dim U_{C_j}(\Gamma,c,\barr{w}) \leq 6d-3$$
unless $e=0$ and $k=1,2$.
Since the linear series determined by the evaluation maps
are base point free, the general intersection
(\ref{erst}) will miss all loci of dimension less than
$6d-2$.

It remains to consider the trees $(\Gamma,c,\barr{w})$ where
$e=0$ and $k=1,2$. If $k=1$, $(\Gamma,c,\barr{w})=\Phi$
satisfies the conditions of Lemma (\ref{aa}).
By Lemma (\ref{aa}),
 $$dim\ I(\Phi,j) \leq 6d-3.$$ Hence, the general intersection
(\ref{erst}) will miss all the loci $U_C(\Phi,c,(0,d))$.

If $k=2$, there are two cases to consider. If there exists
a vertex of weight $d$, then $(\Gamma,c,\barr{w})=\Phi$
satisfies the conditions of Lemma (\ref{bb}).
By Lemma (\ref{bb}),
$$dim\ I(\Phi,j) \leq 6d-4.$$
If $w_1+w_2=d$ is a positive partition, then the
image of every map $[\mu]\in U_C(\Gamma,c,\barr{w})$ is
the union of two rational curves of degrees $w_1$ and $w_2$.
No such unions pass through $3d-1$ general points.
The proof of claim (\ref{erst}) is complete.

For $\infty>j\neq 0,1728$, $\pi_W^{-1}(j)$ is a
reduced, irreducible divisor of $W_{1}(d)$. Since
the linear series determined by the evaluation maps are
base point free, the general intersection cycle
\begin{equation}
\label{al}
Z_{\barr{x}} \cap \pi_W^{-1}(j)
\end{equation}
is a reduced collection of $Z\cap T$ points. The general intersection
cycle $(\ref{al})$ also consists exactly of the reduced, nodal, degree $d$
elliptic plane curves with $j$-invariant $j$ passing through the
points $\barr{x}$.

The argument for $j=0,1728$ is identical except that
$\pi_W^{-1}(0)$, and $\pi_W^{-1}(1728)$ are divisors in $W_1(d)$
with multiplicity $3$, $2$ respectively.
These multiplicities arise from the
extra automorphisms for $j=0,1728$.
Therefore the cycle (\ref{al}) is a collection of
$${1\over 3}\cdot Z\cap T,$$ $${1\over 2}\cdot Z\cap T$$ triple and double
points respectively.
\end{pf}
It remains to evaluate $Z\cap T$.
\begin{lm}
$Z\cap T= {d-1\choose 2} N_d$.
\end{lm}
\begin{pf}

It is first established for a general representative $Z_{\barr{x}}$,
\begin{equation}
\label{lrst}
Z_{\barr{x}} \cap \pi_{\barr{W}}^{-1}(\infty) \subset \pi_W^{-1}(\infty).
\end{equation}
The statement (\ref{lrst}) is proven by considering
the quasi-projective strata of $\barr{M}_{\infty}(d)$.

By arguments of Lemma (\ref{pal}),
all the loci $U_{\infty}(\Gamma,c,\barr{w})$ where
$(\Gamma,c,\barr{w})$ is not the trivial tree are
avoided in the general intersection (\ref{lrst}). Only the
strata $U_{\infty}(\Lambda,\barr{w})$ remain to be considered.
Let $k+1\geq 2$ be the total number of vertices of $\Lambda$.
By the equations for the dimensions of $U_{\infty}(\Lambda,\barr{w})$,
$$dim \ U_{\infty}(\Lambda, \barr{w}) \leq 6d-2-k\leq 6d-3$$
unless all the circuit vertices have weight zero.
If all circuit vertices have weight zero, $k+1\geq 3$.
Now $$dim \ U_{\infty}(\Lambda, \barr{w}) \leq 6d-k \leq 6d-3$$
unless $k=2$.

Only one stable, marked, weighted, graph with $1$-circuit
$(\Lambda, \barr{w})$
need be considered. Vertices $c_1, c_2$ form a weightless circuit.
Vertex $v_3$ is connected to $c_2$ and $w_3=d$. $(\Lambda,\barr{w})=
\Omega$ satisfies the conditions of Lemma (\ref{cc}). Therefore,
$$dim\ I(\Omega, \infty) \leq 6d-4.$$
Claim (\ref{lrst}) is now proven.

The divisor $\pi_W^{-1}(\infty)$ is
reduced and irreducible in $W_{1}(d)$. As above,
\begin{equation}
\label{tal}
Z_{\barr{x}} \cap \pi_W^{-1}(\infty)
\end{equation}
is a reduced collection of $Z\cap T$ points. The general intersection
cycle $(\ref{tal})$ also consists exactly of degree $d$ {\em maps} of the
$1$-nodal rational curve $C_{\infty}$ passing through $\barr{x}$.
The image of such a map must be
one of the $N_d$ degree $d$, nodal, rational plane curves passing
through $\barr{x}$. The number of distinct birational maps (up to isomorphism)
from
$C_{\infty}$ to a ${d-1 \choose 2}$-nodal plane curve is exactly
${d-1 \choose 2}$. Therefore, $Z\cap T= {d-1\choose 2} N_d$.
\end{pf}

\noindent
Department of Mathematics, University of Chicago

\noindent
rahul@@math.uchicago.edu

\end{document}